\documentclass[a4paper,12pt]{amsart}
\usepackage{amsmath,amssymb,amsfonts,latexsym}
\usepackage{indentfirst}

\newtheorem{theorem}{Theorem}[section]
\newtheorem{lemma}[theorem]{Lemma}
\newtheorem{corollary}[theorem]{Corollary}
\newtheorem{proposition}[theorem]{Proposition}

\newcommand{\z}{\mathbb{Z}}

\newcommand{\fd}{\mathbb{F}}

\begin{document}
\title[Kloosterman Sums with Trace One Arguments]{Codes Associated with $O(3,2^r)$ and Power Moments of Kloosterman Sums with Trace One Arguments}
\author{Dae San Kim}

\address{Department of Mathematics, Sogang University, Seoul
121-742, Korea}

\email{dskim@sogang.ac.kr}
\begin{abstract}
We construct a binary linear code $C(O(3,q))$, associated with the
orthogonal group $O(3,q)$. Here $q$ is a power of two. Then we
obtain a recursive formula for the odd power moments of Kloosterman
sums with trace one arguments in terms of the frequencies of weights
in the codes $C(O(3,q))$ and $C(Sp(2,q))$. This is done via Pless
power moment identity and by utilizing the explicit expressions of
Gauss sums for the orthogonal groups.
\end{abstract}

\thanks{This work was supported by National Research Foundation of Korea Grant funded by
the Korean Government 2009-0072514.}

\maketitle
\section {Introduction}

Let $\psi$ be a nontrivial additive character of the finite field
$\fd{_q}$ with $q=p^{r}$ elements ($p$ a prime). Then the
Kloosterman sum $K(\psi;a)$ (\cite{RH}) is defined by

\begin{equation*}
 K(\psi;a)=\sum_{\alpha \in \fd_{q}^{*}}\psi(\alpha+a \alpha^{-1})\\
(a \in \fd_{q}^{*}).
\end{equation*}

The Kloosterman sum was introduced in 1926(\cite{HDK}) to give an
estimate for the Fourier coefficients of modular forms.

For each nonnegative integer $h$, we denote by $MK(\psi)^{h}$ the
$h$-th moment of the Kloosterman sum $K(\psi;a)$, i.e.,

\begin{equation*}
MK(\psi)^{h}=\sum_{a  \in \fd_{q}^{*}} K(\psi ;a)^{h}.
\end{equation*}

If $\psi=\lambda$ is the canonical additive character of $\fd_{q}$,
then $MK(\lambda)^{h}$ will be simply denoted by $MK^{h}$.

Also, we introduce two incomplete power moments of Kloosterman sums,
namely the one with the sum over all $a$ in $\fd_{q}^{*}$ with
$tra=0$ and the other with the sum over all $a$ in $\fd_{q}^{*}$
with $tra=1$. For every nonnegative integer $h$, and $\psi$ as
before, we define

\begin{equation}\label{a1}
T_{0}K(\psi)^{h}=\sum_{a  \in \fd_{q}^{*},\;\;tr a=0}^{}K(
\psi;a)^{h},\;\;T_{1}K(\psi)^{h}= \sum_{a  \in \fd_{q}^{*},\;\;tr
a=1}^{}K(\psi ;a)^{h},
\end{equation}
which will be respectively  called the $h$-th moment of Kloosterman
sums with ``trace zero arguments" and that with ``trace one
arguments." Then, clearly we have
\begin{equation}\label{a2}
MK(\psi)^h = T_{0}K(\psi)^h+T_{1}K(\psi)^h.
\end{equation}

If $\psi=\lambda$ is the canonical additive character of $\fd_{q}$,
then $T_{0}K(\lambda)^{h}$ and $T_{1}K(\lambda)^{h}$ will be
respectively denoted by $T_{0}K^{h}$ and $T_{1}K^{h}$, for brevity.

Explicit computations on power moments of Kloosterman sums were
initiated in the paper \cite{HS} of Sali\'{e} in 1931, where it is
shown that for any odd prime $q$,

\begin{equation*}
MK^h=q^{2}M_{h-1}-(q-1)^{h-1}+2(-1)^{h-1}\;\;\;\;(h \geq 1).
\end{equation*}
Here $M_{0}=0$, and, for $h \in \z_{>0}$,

\begin{equation*}
M_{h}=| \{ (\alpha_{1},\cdots, \alpha_{h}) \in (\fd_{q}^{*})^{h} |
\sum_{j=1}^{h}\alpha_{j} =1=\sum_{j=1}^{h}\alpha_j^{-1}\}|.
\end{equation*}

For $q=p$ an odd prime, Sali\'{e} obtained $MK^1$, $MK^2$, $MK^3$,
$MK^4$ in \cite{HS} by determining $MK^1$, $MK^2$, $MK^3$. $MK^5$
can be expressed in terms of the $p$-th eigenvalue for a weight 3
newform on $\Gamma_{0}(15)$ (cf. \cite{RL}, \cite{CJM}]).  $MK^6$
can be expressed in terms of the $p$-th eigenvalue for a weight 4
newform on $\Gamma_{0}(6)$ (cf. \cite{KJB}). Also, based on
numerical evidence, in \cite{RJ} Evans was led to propose a
conjecture which expresses  $MK^7$ in terms of Hecke eigenvalues for
a weight 3 newform on $\Gamma_{0}(525)$ with quartic nebentypus of
conductor 105.

From now on, let us assume that $q=2^r$. Carlitz\cite{L1} evaluated
$MK^h$ for $h \leq 4$. Recently, Moisio was able to find explicit
expressions of $MK^h$, for $h \leq 10$ (cf. \cite{M1}). This was
done, via Pless power moment identity, by connecting moments of
Kloosterman sums and the frequencies of weights in the binary
Zetterberg code of length $q+1$, which were known by the work of
Schoof and Vlugt in \cite{RS}.

In this paper, we will show the main Theorem \ref{A} giving a
recursive formula for the odd power moments of Kloosterman sums with
trace one arguments. To do that, we construct a binary linear code
$C(O(3,q))$, associated with the orthogonal group $O(3,q)$, and
express those power moments in terms of the frequencies of weights
in the codes $C(O(3,q))$ and $C(Sp(2,q))$. Here $C(Sp(2,q))$ is a
binary linear code associated with the symplectic group $Sp(2,q)$
(cf. \cite{D4}). Then, thanks to our previous results on the
explicit expressions of ``Gauss sums" for the orthogonal group
$O(2n+1,q)$ \cite{DY}, we can express the weight of each codeword in
the dual of the code in terms of Kloosterman sums. Then our formula
will follow immediately from the Pless power moment identity. The
reader is referred to \cite{D2} for the construction  of code
associated with other type of orthogonal group and the derivation of
recursive formulas generating power moments of Kloosterman sums.

In \cite{D4}, for both $n$, $q$ powers of two, a binary linear code
$C(SL(n,q))$ associated with the finite special linear group
$SL(n,q)$ was constructed in order to produce a recursive formula
for the power moments of multi-dimensional Kloosterman sums in terms
of the frequencies of weights in that code. On the other hand, in
\cite{D3}, for $q$ a power of three, two infinite families of
ternary linear codes associated with double cosets in the symplectic
group $Sp(2n,q)$ were constructed in order to generate infinite
families of recursive formulas for the power moments of Kloosterman
sums with square arguments and for the even power moments of those
in terms of the frequencies of weights in those codes.

Henceforth, we agree that the binomial coefficient $\binom{b}{a}=0$,
if $a > b$ or $a < 0$.

\begin{theorem}\label{A}
Let $q=2^r$. Then we have the following. For $h=1,3,5,\cdots,$

\begin{align}\label{a3}
  \begin{split}
 T_{1}K^{h}&=- \sum_{1 \leq j \leq {h-2},\;\;j\;\;odd}^{}\binom{h}{j}(q^{2}-1)^{h-j}T_{1}K^{j}\\
           &+q^{1-h}\sum_{j=0}^{min \{N ,h \}}(-1)^{j}D_j \sum_{t=j}^{h}t!S(h,t)2^{h-t-1} \binom{N
           -j}{N-t},
 \end{split}
\end{align}
where $N=|O(3,q)|=q(q^2 -1)$, and $D_j=C_j-\hat{C_j}$ $(0 \leq j
\leq N)$, with $\{C_{j}\}_{j=0}^{N}$, $\{\hat{C_{j}}\}_{j=0}^{N}$
respectively the weight distributions of $C(O(3,q))$ and
$C(Sp(2,q))$ given by: for $j=0, \cdots,N$,

\begin{equation}\label{a4}
C_j=\sum \binom{q^{2}}{\nu_{1}}\prod_{tr((\beta-1)^{-1}
)=0}\binom{q^{2}+q}{\nu_{\beta }}\prod_{tr((\beta-1)^{-1})=1}
\binom{q^{2}-q}{\nu_{\beta}},
\end{equation}

\begin{equation}\label{a5}
\hat{C_j}=\sum_{}^{}\binom{q^{2}}{\nu_{0}} \prod_{tr(\beta^{-1})=0}
\binom{q^{2}+q}{\nu_{\beta}} \prod_{tr( \beta^{-1})=1}
\binom{q^{2}-q}{\nu_{\beta }}.
\end{equation}
Here the first sum in (\ref{a3}) is 0 if $h=1$ and the unspecified
sums in (\ref{a4}) and (\ref{a5}) run over all the sets of
nonnegative integers $\{\nu_{\beta}\}_{\beta  \in \fd_{q}}$
satisfying $\sum_{\beta \in \fd_{q}} \nu_{\beta}=j$ and $\sum_{\beta
\in \fd_{q}} \nu_{\beta}\beta=0$. In addition, $S(h,t)$ is the
Stirling number of the second kind defined by
\begin{equation}\label{a6}
S(h,t)=\frac{1}{t!}\sum_{j=0}^{t}(-1)^{t-j}\binom{t}{j}j^{h}.
\end{equation}
\end{theorem}
\section{$O(2n+1,q)$}
For more details about this section, one is referred to the paper
\cite{DY}. Throughout this paper, the following notations will be
used:
\begin{itemize}
 \item [] $q = 2^r$ ($r \in \mathbb{Z}_{>0}$),\\
 \item [] $\mathbb{F}_{q}$ = the finite field with $q$ elements,\\
 \item [] $Tr A$ = the trace of $A$ for a square matrix $A$,\\
 \item [] $^tB$ = the transpose of $B$ for any matrix $B$.
\end{itemize}\

Let $\theta$ be the nondegenerate quadratic form on the vector space
$\fd_q^{(2n+1) \times 1}$ of all $(2n+1) \times 1$ column vectors
over $\fd_q$, given by

\begin{equation*}
\theta(\sum_{i=1}^{2n+1} x_i e^i) = \sum_{i=1}^{n}
x_{i}x_{n+i}+x_{2n+1}^{2},
\end{equation*}
where $\{e^1=^t[10\ldots0], e^2=^t[010\ldots
0],\ldots,e^{2n+1}=^t[0\ldots01]\}$ is the standard basis of
$\fd_q^{(2n+1) \times 1}$.

The group $ O(2n+1,q)$ of all isometries of $(\fd_q^{(2n+1) \times
1}, \; \theta)$ consists of the matrices

\[
\left[%
\begin{matrix}
  A & B & 0 \\
  C & D & 0 \\
  g & h & 1\\
\end{matrix}%
\right]%
\;\;\;(A,B,C,D \;\;n \times n, g, h \;\;1 \times n)\\
\]
$~$\\
in $GL(2n+1,q)$ satisfying the relations:

\begin{align*}
  \begin{split}
  &^{t}AC+{}^{t}gg \;\;\; \textmd{is alternating},\\
  &^{t}BD+{}^{t}hh \;\;\; \textmd{is alternating},\\
  &^{t}AD+{}^{t}CB=1_{n}.
  \end{split}
\end{align*}
Here an $n \times n$ matrix $(a_{ij})$ is called alternating if

\begin{equation*}
\begin{cases}
 a_{ii}=0,              & \text{for $1 \leq i \leq n$},\\
 a_{ij}= -a_{ji}=a_{ji}, & \text{for $1 \leq i < j \leq n$.}
\end{cases}
\end{equation*}

Also, one observes, for example, that ${}^{t}AC+{}^{t}gg$  is
alternating if and only if ${}^{t}AC=$ ${}^{t}CA$ and
$g=\sqrt{diag({}^{t}AC)}$, where $\sqrt{diag({}^{t}AC)}$ indicates
the $1 \times n$ matrix $[\alpha_1,\cdots, \alpha_n]$ if the
diagonal entries of ${}^{t}AC$ are given by

\begin{equation*}
(^{t}AC)_{11}=\alpha_1^2,\cdots,\;\;({}^t
AC)_{nn}=\alpha_{n}^{2},\;\; \textmd{for} \;\;\alpha_i \in \fd_q.\\
\end{equation*}
$~$\\
$~$
As is well known, there is an isomorphism of groups

\begin{equation}\label{a7}
\iota: O(2n+1,q)\rightarrow Sp(2n,q)\;\;\;(\left[%
\begin{matrix}
  A & B & 0 \\
  C & D & 0 \\
  g & h & 1\\
\end{matrix}%
\right] \mapsto \left[%
\begin{matrix}
  A & B \\
  C & D \\
\end{matrix}%
\right]).
\end{equation}

In particular, for any $w \in O(2n+1,q)$,
\begin{equation}\label{a8}
Tr w=Tr \iota(w)+1.
\end{equation}
Here the symplectic group  $Sp(2n,q)$ over the field $\fd_{q}$ is
defined as:
\begin{equation*}
Sp(2n,q)=\{w \in GL(2n,q) |^t wJw=J \},
\end{equation*}
with

\begin{equation*}
J=
\left[%
\begin{matrix}
  0 & 1_n \\
  1_n & 0  \\
\end{matrix}%
\right].
\end{equation*}
\begin{equation}\label{a9}
|O(2n+1,q)|=|Sp(2n,q)|=q^{n^2} \prod_{j=1}^n(q^{2j}-1)\;\;\;(cf.\;\;
(7), [4]).
\end{equation}

For integers $n,r$ with $0 \leq r \leq n$, the $q$-binomial
coefficients are defined as:

\begin{equation}\label{a10}
\left[ \substack{n \\ r}
 \right]_q = \prod_{j=0}^{r-1} (q^{n-j} - 1)/(q^{r-j}-1).
 \end{equation}

\section{Gauss sums for $O(2n+1,q)$}
The following notations will be used throughout this paper.
\begin{gather*}
tr(x)=x+x^2+\cdots+x^{2^{r-1}} \text{the trace function}
~\mathbb{F}_{q}
\rightarrow \mathbb{F}_2,\\
\lambda(x) = (-1)^{tr(x)} ~\text{the canonical additive character
of} ~\fd_q.
\end{gather*}
Then any nontrivial additive character $\psi$ of $\fd_q$ is given by
$\psi(x) = \lambda(ax)$, for a unique $a \in \fd_q^*$.\\

For any nontrivial additive character $\psi$ of $\fd_q$ and $a \in
\fd_q^*$, the Kloosterman sum $K_{GL(t,q)}(\psi ; a)$ for $GL(t,q)$
is defined as
\begin{equation}\label{a11}
K_{GL(t,q)}(\psi ; a) = \sum_{w \in GL(t,q)} \psi(Trw + a~Trw^{-1}).
\end{equation}
Observe that, for $t=1 $, $ K_{GL(1,q)}( \psi;a)$ denotes the
Kloosterman sum $K(\psi;a) $.

In \cite{D1}, it is shown that $K_{GL(t,q)}(\psi ; a)$ ~satisfies
the following recursive relation: for integers $t \geq 2$, ~$a \in
\fd_q^*$ ,
\begin{multline}\label{a12}
K_{GL(t,q)}(\psi ; a) = q^{t-1}K_{GL(t-1,q)}(\psi ; a)K(\psi
;a)\\
+ q^{2t-2}(q^{t-1}-1)K_{GL(t-2,q)}(\psi ; a),
\end{multline}
where we understand that $K_{GL(0,q)}(\psi ; a)=1$.\\

From \cite{D1} and \cite{DY},  the Gauss sum for $O(2n+1,q)$ is
equal to $\psi(1)$ times that for $Sp(2n,q)$ and is  given by (cf.
(\ref{a10}), (\ref{a11})):

\begin{align*}
  \begin{split}
 &\sum_{w \in O(2n+1,q)} \psi(Tr w)=\psi(1) \sum_{ w \in Sp(2n,q)}\psi(Trw)\\
 &=\psi(1)q^{\binom{n+1}{2}}\sum_{0\leq r\leq n,\;\; r \;\;\textmd{even}} q^{rn- \frac{1}{4 }r^2}\left[ \substack{n \\ r}
 \right]_q  \prod_{j=1}^{r/2} (q^{2j-1} -1)K_{GL(n-r,q)}(\psi;1).
  \end{split}
\end{align*}

Here $\psi$ is any nontrivial additive character of $\fd_{q}$. For
our purposes, we only need the following expression of the Gauss sum
for $O(3,q)$. So we state it separately as a theorem.

\begin{theorem}\label{B}
Let $\psi$ be any nontrivial additive character of $\fd_{q}$. Then
we have
\begin{equation*}
\sum_{w \in O(3,q)} \psi(Tr w)=\psi(1)qK(\psi;1).
\end{equation*}
\end{theorem}

\begin{proposition}(\cite{D2})\label{C}
For $n=2^s$ ($s \in \z_{\geq 0}$), and $\lambda$ the canonical
additive character of $\fd_q$,

\begin{equation*}
K( \lambda ;a^n )=K(\lambda;a).
\end{equation*}
\end{proposition}

The next corollary follows from Theorem 2 and Proposition 3 and by
simple change of variables.

\begin{corollary}\label{D}
 Let $\lambda$ be the canonical additive character of $\fd_q$, and let $a \in \fd_{q}^{*}$. Then we have
\begin{equation}\label{a13}
\sum_{ w \in O(3,q)} \lambda(aTr w)=\lambda(a)qK(\lambda ;a ).
\end{equation}
\end{corollary}
\begin{proposition}\label{E}(\cite{D2})
Let $\lambda $ be the canonical additive character of $\fd_{q}$,
$\beta \in \fd_{q}$. Then

\begin{align}\label{a14}
  \begin{split}
 &\sum_{a \in \fd_{q}^{*}} \lambda(-a \beta)K(\lambda ;a)\\
 &=\begin{cases}
    q\lambda (\beta^{-1})+1, & \hbox{if $\beta \neq 0$,} \\
      1, & \hbox{if $\beta = 0$.} \\
\end{cases}
\end{split}
\end{align}
\end{proposition}

Let $G(q)$ be  one of finite classical groups over $\fd_{q}$. Then
we put, for each $\beta \in \fd_{q}$,

\begin{equation*}
N_{G(q)}(\beta)=|\{w \in G(q) |Tr(w)= \beta \} |.
\end{equation*}

Then it is easy to see that

\begin{equation}\label{a15}
qN_{G(q)}(\beta)=|G(q)|+ \sum_{a \in \fd_q^*}\lambda(-a
\beta)\sum_{w \in G(q)}\lambda(aTr w).
\end{equation}

For brevity, we write

\begin{equation}\label{a16}
n(\beta)=N_{O(3,q)}(\beta).
\end{equation}

Using (\ref{a13})-(\ref{a15}), one derives the following.

\begin{proposition}\label{F}
We have
\begin{equation}\label{a17}
 n(\beta)
 =\begin{cases}
    q^2,    & \textmd{if}\;\; \beta =1 ,\\
     q^2+q , &\textmd{if} \;\; tr((\beta-1)^{-1})=0, \\
     q^2-q , &\textmd{if} \;\; tr((\beta-1)^{-1})=1.
     \end{cases}
\end{equation}
\end{proposition}

\begin{corollary}\label{G}
$Tr : O(3,q) \rightarrow \fd_{q}$ is surjective.
\end{corollary}
\proof From (\ref{a17}), we see that $n(\beta) > 0$, for all $\beta
\in \fd_{q}$. Alternatively, this also follows from the easily shown
surjectivity of $Tr:Sp(2n,q) \rightarrow \fd_q$ and the facts in
(\ref{a7}) and (\ref{a8}).\;\;\;$\square$\\

In below, ``the sum over $tr a=0$" will mean ``the sum over all
nonzero $a \in \fd_{q}^{*}$, with $tr a=0$."

The result (b) in below and hence (a) follow also from (\ref{a3}).

\begin{proposition}\label{H}
Let $q=2^r$. Then we have the following.\\

(a) $T_{0}K=1+\frac{1}{2}(-1)^{r}q$.\\

(b) $T_{1}K=\frac{1}{2}(-1)^{r+1}q$.
\end{proposition}

\proof (a)
\begin{align*}
  \begin{split}
 T_{0}K&=\sum_{tra=0}\sum_{x \in \fd_q^*}\lambda(x+ax^{-1})\\
       &=\sum_{tr a=0}\sum_{x \in \fd_q^*}\lambda(x^{-1}+ax)\\
       &=\frac{1}{2}\sum_{x \in \fd_q^*}\lambda(x^{-1})\sum_{\alpha \in \fd_q \backslash \{0,1\}}\lambda((\alpha^2+\alpha)x)\\
       &=\frac{1}{2}\sum_{x \in \fd_q^*}\lambda(x^{-1})\{\sum_{\alpha \in \fd_q}\lambda(x \alpha^2+x \alpha)-2 \}
  \end{split}
\end{align*}
\begin{equation}\label{a18}
=1+\frac{1}{2}\sum_{x \in \fd_q^*}\lambda(x^{-1})\sum_{\alpha \in
\fd_q}\lambda(x \alpha^2+x \alpha).
\end{equation}

Now, from Theorem 5.34 of \cite{RH}, we have

\begin{equation*}
\sum_{\alpha \in \fd_q}\lambda(x \alpha^2+x \alpha)
 =\begin{cases}
   q, &\textmd{if}\;\;\; x^2+x=0\;\;(\textmd{i.e},\;\; x \in \fd_2),\\
    0 ,& \textmd{otherwise}.
\end{cases}
\end{equation*}

Thus (\ref{a18}) equals
\begin{equation*}
1+ \frac{q}{2}\sum_{x \in
\fd_2^{*}}\lambda(x^{-1})=1+\frac{1}{2}(-1)^{r}q.
\end{equation*}
(b) This follows from (a), since $T_{1}K=MK-T_{0}K=1-T_{0}K$.\;\;\;$
\square$

\section{Construction of codes}
Let
\begin{equation}\label{a19}
N=|O(3,q)|=q(q^2 -1).
\end{equation}
Here we will construct a binary linear code $C(O(3,q))$ of length
$N$, associated with the orthogonal group $O(3,q)$.

Let $g_{1}, g_{2}, \cdots, g_{N}$ be a fixed ordering of the
elements in the group $O(3,q)$. Also, we put
\begin{equation*}
v=(Tr g_{1},Tr g_{2},\cdots,Tr g_{N}) \in \fd_q^N.
\end{equation*}
Then the binary linear code $C(O(3,q))$ is defined as

\begin{equation}\label{a20}
C(O(3,q))=\{ u \in \fd_{2}^N  |u \cdot v =0 \},
\end{equation}
where the dot denotes the usual inner product in $\fd_{q}^{N}$.

The following Delsarte's theorem is well-known.

\begin{theorem}(\cite{FN})\label{I}
Let $B$ be a linear code over $ \fd_{q}$. Then
\begin{equation*}
(B|_{\fd_{2}})^{\bot }=tr(B^{\bot}).
\end{equation*}
\end{theorem}

In view of this theorem, the dual $C(O(3,q))^{\bot}$ is given by
\begin{equation}\label{a21}
C(O(3,q))^{\bot}= \{c(a)=(tr(aTrg_{1}), \cdots,tr(aTrg_{N}))|a \in
\fd_q \}.
\end{equation}

\begin{proposition}\label{J}
For every $q=2^r$, the map $\fd_{q} \rightarrow C(O(3,q))^{\bot}$
($a \mapsto c(a)$) is an $\fd_2$-linear isomorphism.
\end{proposition}

\proof The map is clearly $\fd_2$-linear and surjective. Let $a$ be
in the kernel of the map. Then, in view of Corollary \ref{G},
$tr(\alpha\beta)=0$, for all $\beta \in \fd_q$. Since the trace
function $\fd_q \rightarrow \fd_2$ is surjective,
$a=0$.\;\;\;$\square$

\section{Power moments of Kloosterman sums with trace one arguments}
In this section, we will be able to find, via Pless power moment
identity, a recursive formula for the power moments of Kloosterman
sums with trace one arguments in terms of the frequencies of weights
in the codes $C(O(3,q))$ and $C(Sp(2,q))$.

\begin{theorem}\label{K}(Pless power moment identity)
Let $ B$ be an $q$-ary $[n,k]$ code, and let $B_{i}$(resp.$B_{i}
^{\bot})$ denote the number of codewords of weight $i$ in $B$(resp.
in $B^{\bot})$. Then, for $h=0,1,2, \cdots$,

\begin{align}\label{a22}
\sum_{j=0}^{n}j^{h}B_{j}=\sum_{j=0}^{min \{ n,h \}}(-1)^{j}B_{j}
^{\bot} \sum_{t=j}^{h} t! S(h,t)q^{k-t}(q-1)^{t-j}\binom{n-j}{n-t},
\end{align}
where $S(h,t)$ is the Stirling number of the second kind defined in
(\ref{a6}).
\end{theorem}

\begin{lemma}\label{L}
Let $c(a)=(tr(aTrg_{1}), \cdots,tr(aTrg_{N})) \in C(O(3,q))^{\bot}$,
for $a \in \fd_q^{*}$. Then the Hamming weight $w(c(a))$ can be
expressed as follows:
\begin{equation}\label{a23}
w(c(a))=\frac{1}{2}q \{(q^2-1)-\lambda(a)K(\lambda;a)\}.
\end{equation}
\end{lemma}
\proof
\begin{align*}
  \begin{split}
  w(c(a))&=\frac{1}{2}\sum_{j=1}^{N}(1-(-1)^{tr(aTrg_j)})\\
             &=\frac{1}{2}(N- \sum_{w \in O(3,q)}\lambda(a Tr w)).
  \end{split}
\end{align*}
Our results now follow from (\ref{a19}) and (\ref{a13}).
\;\;\;$\square$

Let $u=(u_1,\cdots,u_N) \in \fd_2^N$, with $\nu_\beta$ 1's in the
coordinate places where $ Tr(g_j)=\beta$, for each $\beta \in
\fd_q$. Then we see from the definition of the code $C(O(3,q))$ (cf.
(\ref{a20})) that $u$ is a codeword with weight $j$ if and only if
$\sum _{\beta \in \fd_{q}} \nu_{\beta} =j$ and $\sum_{\beta \in
\fd_{q}}\nu _{\beta}\beta =0$ (an identity in $\fd_{q}$). Note that
there are $\prod_{\beta \in \fd_q}\binom{n(\beta)}{\nu_\beta}$ many
such codewords with weight $j$. Now, we obtain the following formula
in (\ref{a24}), by using the explicit value of $n(\beta)$ in
(\ref{a17}).

\begin{theorem}\label{M}
Let $q=2^r$ be as before, and let $\{C_j \}_{j=0}^N$ be the weight
distribution of $C(O(3,q))$. Then, for $j=0, \cdots,N$,

\begin{equation}\label{a24}
C_j=\sum \binom{q^{2}}{\nu_{1}} \prod_{tr((\beta-1)^{-1})=0}\binom{q
^{2}+q}{\nu_{\beta}}\prod_{tr((\beta-1)^{-1})=1} \binom{q^{2}-q}{\nu
_{\beta}},
\end{equation}
where the sum runs over all the sets of nonnegative integers $\{\nu
_{\beta}\}_{\beta  \in \fd_{q}}$ satisfying

\begin{equation}\label{a25}
\sum_{\beta \in \fd_{q}} \nu_{\beta }=j \;\;\; \textmd{and}\;\;\;
\sum_{\beta \in \fd_{q}} \nu_{\beta }\beta =0,
\end{equation}
and the first and second products are over all $\beta \neq 1$,
respectively with $tr((\beta-1)^{-1})=0$ and $tr(
(\beta-1)^{-1})=1$.
\end{theorem}

\begin{corollary}\label{N}
Let $\{C_j \}_{j=0}^N$ be the weight distribution of $C(O(3,q)$.
Then we have $C_j=C_{N-j}$, for all $j$, with $0 \leq j \leq N $.
\end{corollary}
\proof Under the replacements $\nu_{\beta} \rightarrow n(\beta)-\nu_
{\beta}$, for each $\beta \in \fd_q$, the first equation in
(\ref{a25}) is changed to $N-j$, while the second one in (\ref{a25})
and the summands in (\ref{a24}) are left unchanged. Here the second
sum in (\ref{a25}) is left unchanged, since $\sum_{\beta \in \fd_q }
n(\beta)\beta=0$, as one can see by using the explicit expression of
$n(\beta)$ in (\ref{a17}). \;\;\;$\square$\\

The recursive formula in the following theorem follows from the
study of codes associated with symplectic group $Sp(2,q)=SL(2,q)$.
They are slightly modified from their original versions, which makes
them more usable in below.

\begin{theorem}\label{O}(\cite{D4})
For $h=1,2,3, \cdots$,

\begin{align}\label{a26}
  \begin{split}
  (\frac{q}{2})^h  \sum_{j=0}^{ h}&(-1)^j \binom{h}{j}(q^2-1)^{h-j}MK^j\\
   &=q\sum_{j=0}^{min \{N ,h \}}(-1)^{j}\hat{C_{j}} \sum_{t=j}^{h} t!S(h,t)2^{-t}\binom{N -j}{N-t},
  \end{split}
\end{align}
where $N=q(q^{2}-1)=|Sp(2,q)|=|O(3,q)|$, $S(h,t)$ indicates the
Stirling number of the second kind as in (\ref{a6}), and
$\{\hat{C_{j}}\}_{j=0}^{N}$ denotes the weight distribution of the
code $C(Sp(2,q))$, given by
\begin{equation*}
\hat{C_{j}}=\sum \binom{q^{2}}{\nu_{0}} \prod_{tr(\beta^{-1}
)=0}\binom{q^{2}+q}{\nu_{\beta}}\prod_{tr(\beta^{-1})=1} \binom{q
^{2} -q}{\nu_{\beta }}\;\;\;(0 \leq j \leq N).
\end{equation*}
Here the sum runs over all the sets of nonnegative integers $\{\nu
_{\beta}\}_{\beta  \in \fd_{q}}$ satisfying $\sum_{\beta \in
\fd_{q}} \nu_{\beta}=j$ and $\sum_{\beta  \in \fd_{q}} \nu_{\beta}
\beta =0$, and the first and second product run over the elements
$\beta \in \fd_{q}^{*}$, respectively with $tr(\beta^{-1})=0$ and
$tr(\beta^{-1} )=1$.

\end{theorem}

We are now ready to apply the Pless power moment identity in
(\ref{a22}) to $C(O(3,q))^{\bot}$, in order to obtain the result in
Theorem \ref{A} (cf. (\ref{a3})) about a recursive formula.

Then the left hand side of that identity in (\ref{a22}) is equal to

\begin{equation}\label{a27}
\sum_{ a \in F_q^{*} } w(c(a))^h,
\end{equation}
with the $w(c(a))$ given by (\ref{a23}).

(\ref{a27}) is now given by

\begin{align*}
  \begin{split}
 &(\frac{q}{2})^h \sum_{a \in \fd_q^*}(q^2-1-\lambda (a)K(\lambda ;a))^h\\
 &=(\frac{q}{2})^h \sum_{tr a=0}(q^2-1-K (\lambda ;a))^h+(\frac{q}{2})^h \sum_{tr a=1}(q^2-1+K ( \lambda ;a))^h\\
 &=(\frac{q}{2 })^h \sum_{tr a=0} \sum_{j=0}^{ h}(-1)^{j} \binom{h}{j}(q^2-1)^{h-j}(K( \lambda ;a))^j\\
 &\qquad \qquad +(\frac{q}{2})^h  \sum_{tr a=1} \sum_{j=0}^{h}\binom{h}{j}(q^2-1)^{h-j}(K(\lambda;a))^j\\
 &=(\frac{q}{2})^h \sum_{j=0}^{h}(-1)^{j}\binom{h}{j}(q^2-1)^{h-j}(MK^j-T_1 K^j)\;\;\;(\psi= \lambda \;\;\textmd{case of}\;\; (1),(2))\\
 &\qquad \qquad +(\frac{q}{2})^h \sum_{j=0}^{ h}\binom{h}{j}(q^2-1)^{h-j}T_1K^j\\
 &=(\frac{q}{2})^h \sum_{j=0}^{ h}(-1)^j \binom{h}{j}(q^2-1)^{h-j}MK^j\\
 & \qquad \qquad +2(\frac{q}{2})^h \sum_{0 \leq j \leq h,\;\;j\;\;\textmd{odd }} \binom{h}{j}(q^2-1)^{h-j}T_1 K^j
  \end{split}
\end{align*}
\begin{align}\label{a28}
  \begin{split}
 &=q \sum_{j=0}^{min \{N ,h \}}(-1)^{j} \hat{C_{j}}\sum_{t=j}^{h}t! S(h,t)2^{-t}\binom{N -j}{N-t} \;\;\;(\textmd{cf}.\;\;(26)) \qquad\qquad\\
 &\qquad \qquad  +2 (\frac{q}{2})^h \sum_{0 \leq j \leq h,\;\;j \;\;\textmd{odd}} \binom{h}{j}(q^2-1)^{h-j} T_1 K^j.
  \end{split}
\end{align}

Here one has to separate the term corresponding to $l=h$ in
(\ref{a28}), and note $dim_{\fd_2}C(O(3,q))=r$.


\end{document}